\documentclass[11pt]{article}
\usepackage[latin1]{inputenc}
\usepackage{amssymb}
\usepackage{amsmath}
\usepackage{amsfonts}
\usepackage{amsthm}
\usepackage{enumerate}
\usepackage{bbm}

\title{A note on the Cauchy problem for the 2D generalized Zakharov-Kuznetsov equations}
\date{}
\author{Francis Ribaud\thanks{Laboratoire d'Analyse et de Math\'ematiques Appliqu\'ees - Universit\'e Paris-Est - 5 Bd. Descartes - Champs-Sur-Marne, 77454 Marne-La-Vall\'ee Cedex 2 (FRANCE)} \and St\'ephane Vento\thanks{Laboratoire Analyse, G\'eom\'etrie et Applications - Universit\'e Paris 13 - Institut Galil\'ee, 99 avenue J.B. Cl\'ement,93430 Villetaneuse (FRANCE)}}

\numberwithin{equation}{section}

\newtheorem{theorem}{Theorem}[section]

\newcommand\re[1]{(\ref{#1})}
\def\R{\mathbb{R}}

\newcommand\cro[1]{\langle #1 \rangle}
\begin{document}

\maketitle
\noindent {\bf Abstract.}\, In this note we study the generalized 2D Zakharov-Kuznetsov equations
$\partial_tu+\Delta\partial_xu+u^k\partial_xu=0$ for $k\ge 2$. By an iterative method we prove the local well-posedness of these equations in the
Sobolev spaces $H^s(\mathbb{R}^2)$ for $s>1/4$ if $k=2$, $s>5/12$ if $k=3$ and $s>1-2/k$ if $k\ge 4$.
\\

\noindent
{\bf Keywords:} KdV-like equations, Cauchy problem\\
{\bf AMS Classification:} 35Q53, 35B65, 35Q60

\section{Introduction}
In this short note, we are interested with the Cauchy problem associated to the generalized Zakharov-Kuznetsov (gZK) equations
\begin{equation}\label{gzk}
u_t+\Delta u_x+u^ku_x=0,
\end{equation}
in two-dimensional space and for $k=2,3,4,...$. These equations are natural multi-dimensional generalizations of the well-known generalized Korteweg-de Vries equations and have been derived in \cite{ZK} when $k=1$ to model the propagation of nonlinear ionic-sonic waves in a magnetized plasma.

We give sharp results concerning the well-posedness issue in standard Sobolev spaces $H^s(\R^2)$ for suitable $s\in\R$.
This work follows and use similar technics as in  the paper \cite{RV} where we proved that the 3D associated problem for $k=1$ is locally well-posed in $H^{1^+}(\R^3)$.

Remark that the Sobolev spaces $\dot{H}^s(\R^2)$ are invariant by the natural rescaling of the equation as soon as $s=s_k:=1-2/k$. Thus a natural question is whether (gZK) is well-posed in $H^s(\R^2)$ for $s>s_k$.

\begin{theorem}\label{th}
For any $u_0\in H^s(\R^2)$ with
\[
\left\{\begin{array}{lll}s> 1/4&\textrm{if }k=2,\\ s> 5/12&\textrm{if }k=3,\\ s>1-2/k&\textrm{if }k\ge 4,\end{array}\right.
\]
there exist $T>0$, a Banach space $X^s_T$ and a unique solution $u$ of the Cauchy problem associated to \re{gzk} with $u(0)=u_0$ such that
$u\in X_T^s\cap C_b([0,T], H^s(\R^2)).$
Moreover, the flow-map $u_0\mapsto u$ is Lipschitz on every bounded set of $H^s(\R^2)$.
\end{theorem}
This theorem improves the recent works of Farah, Linares and Pastor in \cite{FLP}-\cite{LP}-\cite{LP2} where local well-posedness was obtained in $H^s(\R^2)$ for $s>3/4$ if $2\le k\le 8$ and $s>s_k$ if $k>8$.

In view of the ill-posedness result obtained in \cite{LP2}, Theorem \ref{th} is optimal (up to the end point) for $k\ge 4$ whereas in the particular cases $k=2$ and 3, we still have a gap (respectively $1/4$ and $1/12$) compared with the scaling index. Concerning the end point $s=s_k$, local well-posedness could perhaps be reached by following the strategy developed in \cite{V}, but with a flow-map only continuous.

In a standard way our proof is based on a fixed point scheme applied to the Duhamel formulation of \re{gzk}:
\begin{equation}\label{duhamel}
u(t) = U(t)u_0 - \frac 1{k+1}\int_0^tU(t-t')\partial_x(u^{k+1})(t')dt'.
\end{equation}
where $U(t)=e^{-t\Delta\partial_x}$ denotes the propagator associated with the linear part of \re{gzk}. Following the works of Kenig, Ponce and Vega \cite{KPV} on the KdV equation,
we use in a crucial way some sharp dispersive estimates for  free solutions. More precisely, these estimates are the well known Kato smoothing effect
\begin{equation}\label{smooth-est}
\|\nabla U(t)\varphi\|_{L^\infty_xL^2_{yT}}\lesssim \|\varphi\|_{L^2},
\end{equation}
which allows to gain one derivative in each spatial direction, and the maximal in time estimate
\begin{equation}\label{maxl4-est}
\|U(t)\varphi\|_{L^4_xL^\infty_{yT}}\lesssim \|\varphi\|_{H^s},\quad s>3/4,
\end{equation}
proved in \cite{LP}. On the other hand, similarly to the generalized KdV equations, the previous bound is no more sufficient to deal with low non-linearities ($k= 2 , 3$)
and we need the following $L^2_x$ based maximal in time estimate (see \cite{F})
\begin{equation} \label{maxl2-est}
\|U(t)\varphi\|_{L^2_xL^\infty_{yT}}\lesssim \|\varphi\|_{H^s},\quad s>3/4.
\end{equation}

\section{Proof of the main result}
\subsection{The case $k\ge 4$}
As mentioned in the introduction, we want to take advantage of the $L^4_xL^\infty_{yT}$ linear estimate \re{maxl4-est}. This motivates the choice of our resolution space:
\[
X^s_T = \{u\in C_b([0,T], H^s(\R^2)): \|u\|_{X^s_T}<\infty\},
\]
where
\[
\|u\|_{X^s_T} = \|u\|_{L^\infty_T H^s_{xy}} + \|\cro{\nabla}^{s+1}u\|_{L^\infty_xL^2_{yT}}+\|\cro{\nabla}^{s-3/4^+}u\|_{L^4_xL^\infty_{yT}}.
\]
Combining estimates \re{smooth-est}-\re{maxl4-est} as well as the straightforward bound
\begin{equation}
\|U(t)\varphi\|_{L^\infty_TL^2}\lesssim \|\varphi\|_{L^2},
\end{equation}
we get
\begin{equation}\label{linhom1-est}
\|U(t)\varphi\|_{X^s_T}\lesssim \|\varphi\|_{H^s}.
\end{equation}
Note that the bound for the second term can be handled by using a low-high frequencies decomposition and next Bernstein inequality and estimate \re{smooth-est}. Having the linear part under control, it remains to deal with the integral term. It is not too hard to adapt the proofs of Propositions 3.5-3.6-3.7 in \cite{RV} to deduce
\begin{equation}\label{linnonhom1-est}
\left\|\int_0^tU(t-t')\partial_xu^{k+1}(t')dt'\right\|_{X^s_T}\lesssim \|\cro{\nabla}^{s-1}\partial_xu^{k+1}\|_{L^1_xL^2_{yT}}.
\end{equation}
The multi-dimensional version of Theorem A.13 in \cite{KPV} applies and leads to
\begin{eqnarray}
\nonumber \|\cro{\nabla}^{s-1}\partial_xu^{k+1}\|_{L^1_xL^2_{yT}} &\lesssim &\|\cro{\nabla}^su^{k+1}\|_{L^1_xL^2_{yT}}\\
\label{liebniz-est} &\lesssim &\|\cro{\nabla}^su\|_{L^{7^+}_xL^{14/3^-}_{yT}}\|u\|_{L^{7k/6^-}_xL^{7k/2^+}_{yT}}^k.
\end{eqnarray}
We claim that the first product in the right hand side of \re{liebniz-est} is controlled by the $X^s_T$ norm of $u$. Indeed, an interpolation argument shows that
\begin{equation}\label{interp-est}
\|\cro{\nabla}^\alpha u\|_{L^p_xL^q_{yT}}\lesssim \|u\|_{X^s_T}
\end{equation}
as soon as there exists $\theta\in [0,1]$ such that
\[
\frac 1p = \frac{1-\theta}4,\quad \frac 1q=\frac{\theta}2,\quad \alpha = \left(s+\frac{7\theta-3}{4}\right)^-.
\]
Taking $\alpha=s$, i.e. $\theta=3/7^+$, it follows that
\[
\|\cro{\nabla}^s u\|_{L^{7^+}_xL^{14/3^-}_{yT}} \lesssim \|u\|_{X^s_T}.
\]
If we choose now $\theta = 4/7k$ in \re{interp-est}, we infer
\[
\|\cro{\nabla}^{(s-3/4+1/k)^-}u\|_{L^{(\frac 14-\frac 1{7k})^{-1}}_xL^{7k/2}_{yT}}\lesssim \|u\|_{X^s_T}.
\]
In order to get the desired contraction factor, we will interpolate this inequality with the bound
\begin{equation}\label{contract-est}
\|\cro{\nabla}^{(s-1)^+}u\|_{L^{\infty^-}_{xyT}}\lesssim T^{0^+}\|\cro{\nabla}^s u\|_{L^\infty_TL^2_{xy}}\lesssim T^{0^+}\|u\|_{X^s_T}.
\end{equation}
This leads to
\begin{equation}\label{nonlin1-est}
\|\cro{\nabla}^{(s-3/4+1/k)^-}u\|_{L^{((\frac 14-\frac 1{7k})^{-1})^+}_xL^{7k/2^+}_{yT}}\lesssim T^{0^+}\|u\|_{X^s_T}.
\end{equation}
By virtue of the Sobolev inequalities, we get
\begin{eqnarray}
\nonumber \|u\|_{L^{7k/6^-}_xL^{7k/2^+}_{yT}} &\lesssim &\|\cro{\nabla}^{(1/4-1/k)^+}u\|_{L^{((\frac 14-\frac 1{7k})^{-1})^+}_xL^{7k/2^+}_{yT}}\\
\nonumber &\lesssim &\|\cro{\nabla}^{(s-3/4+1/k)^-}u\|_{L^{((\frac 14-\frac 1{7k})^{-1})^+}_xL^{7k/2^+}_{yT}}\\
\label{nonlin2-est} &\lesssim & T^{0^+}\|u\|_{X^s_T}
\end{eqnarray}
for $s-3/4+1/k > 1/4-1/k$, that is $s>s_k$. Gathering together \re{linhom1-est}-\re{linnonhom1-est}-\re{liebniz-est}-\re{nonlin1-est} and \re{nonlin2-est} we infer that
\[
\|F(u)\|_{X^s_T}\lesssim \|u_0\|_{H^s}+T^{0^+}\|u\|_{X^s_T}^{k+1},
\]
where $F(u)$ denote the right hand side of \re{duhamel}. The well-posedness result follows then from standard arguments.
\subsection{The case $k=2$}
The proof in this case follows the same lines as in the case $k\ge 4$, but with the $L^4_x$ norm replaced with a $L^2_x$ maximal in time norm. So let us endow the $X^s_T$ space with the norm
\[
\|u\|_{X^s_T} = \|u\|_{L^\infty_T H^s_{xy}} + \|\cro{\nabla}^{s+1}u\|_{L^\infty_xL^2_{yT}}+\|\cro{\nabla}^{s-3/4^+}u\|_{L^2_xL^\infty_{yT}},
\]
for any $s> 1/4$. Using now \re{maxl2-est}, we easily see that
\begin{equation}\label{linhom2-est}
 \|U(t)u_0\|_{X^s_T}\lesssim \|u_0\|_{H^s},
\end{equation}
and
\begin{equation}\label{linnonhom2-est}
 \left\|\int_0^tU(t-t')\partial_xu^3(t')dt'\right\|_{X^s_T}\lesssim \|\cro{\nabla}^{s-1}\partial_xu^3\|_{L^1_xL^2_{yT}}.
\end{equation}
Again, the fractional Leibniz rule yields the bound
\[
\|\cro{\nabla}^{s-1}\partial_xu^3\|_{L^1_xL^2_{yT}} \lesssim \|\cro{\nabla}^su\|_{L^{7/2^+}_xL^{14/3^-}_{yT}} \|u\|_{L^{14/5^-}_xL^{7^+}_{yT}}^2.
\]
By interpolation, we get
\begin{equation}\label{interp2-est}
\|\cro{\nabla}^\alpha u\|_{L^p_xL^q_{yT}}\lesssim \|u\|_{X^s_T}
\end{equation}
for $\alpha$, $p$ and $q$ satisfying
\[
\frac 1p = \frac{1-\theta}2,\quad \frac 1q=\frac{\theta}2,\quad \alpha = \left(s+\frac{7\theta-3}{4}\right)^-
\]
for $0\le\theta\le 1$. On one hand, we deduce that
\[
\|\cro{\nabla}^{s-1}\partial_xu^3\|_{L^1_xL^2_{yT}} \lesssim \|u\|_{X^s _T},
\]
were we took $\theta=3/7^+$ in \re{interp2-est}. On the other hand, for $\theta=2/7^-$, we infer
\[
\|\cro{\nabla}^{(s-1/4)^-}u\|_{L^{14/5^-}_xL^{7^+}_{yT}}\lesssim \|u\|_{X^s_T},
\]
which interpolated with \re{contract-est} gives
\[
\|\cro{\nabla}^{(s-1/4)^-}u\|_{L^{14/5^-}_xL^{7^+}_{yT}}\lesssim T^{0^+}\|u\|_{X^s_T}.
\]
This yields the desired result for $s> 1/4$.

\subsection{The case $k=3$}
Now we consider the intermediate case $k=3$. To prove our result, we define the resolution space as the intersection of the two previous spaces, i.e. equipped with the norm
\[
\|u\|_{X^s_T} = \|u\|_{L^\infty_T H^s_{xy}} + \|\cro{\nabla}^{s+1}u\|_{L^\infty_xL^2_{yT}}+\|\cro{\nabla}^{s-3/4^+}u\|_{(L^2_x\cap L^4_x)L^\infty_{yT}}.
\]
Actually we don't require the full range $L^2_x\cap L^4_x$ for the maximal in time norm, but only
\begin{equation}\label{maxl3-est}
\|\cro{\nabla}^{s-3/4^+}u\|_{L^3_xL^\infty_{yT}}\lesssim \|u\|_{X^s_T}.
\end{equation}
According to \re{linhom1-est}-\re{linnonhom1-est}-\re{linhom2-est}-\re{linnonhom2-est} we have
\[
\|F(u)\|_{X^s_T}\lesssim \|u_0\|_{H^s}+\|\cro{\nabla}^su^4\|_{L^1_xL^2_{yT}}.
\]
Using again the Leibniz rule for fractional derivatives, we infer
\[
\|\cro{\nabla}^su^4\|_{L^1_xL^2_{yT}} \lesssim \|\cro{\nabla}^su\|_{L^{21/4^+}_xL^{14/3^-}_{yT}} \|u\|_{L^{63/17^-}_xL^{21/2^+}_{yT}}^3.
\]
From an interpolation argument with \re{maxl3-est}, we easily check that both these norms are acceptable as soon as $s>5/12$.


\begin{thebibliography}{99}

\bibitem{F} A. V. Faminskii,
{\sl The Cauchy problem for the Zakharov-Kuznetsov equation.} Differential Equations 31 (1995), no. 6, 1002--1012

\bibitem{FLP} L.G. Farah, F. Linares and A. Pastor,
{\sl A note on the 2D generalized Zakharov-Kuznetsov equation: local, global, and scattering results.} 2011, arXiv:1108.3714

\bibitem{KPV} C E. Kenig, G. Ponce and L. Vega,
{\sl Well-posedness and scattering results for the generalized Korteweg-de Vries equation via the contraction principle.} Comm. Pure Appl. Math. 46 (1993), no. 4, 527--620.

\bibitem{LP} F. Linares and A. Pastor,
{\sl Well-posedness for the two-dimensional modified Zakharov-Kuznetsov equation.} SIAM J. Math. Anal. 41 (2009), no. 4, 1323--1339.

\bibitem{LP2} F. Linares and A. Pastor,
{\sl Local and global well-posedness for the 2D generalized Zakharov-Kuznetsov equation.} J. Funct. Anal. 260 (2011), no. 4, 1060--1085.

\bibitem{RV} F. Ribaud and S. Vento,
{\sl Well-posedness results for the 3D Zakharov-Kuznetsov equation.} Preprint, arXiv:1111.2850

\bibitem{V} S. Vento,
{\sl Well-posedness for the generalized Benjamin-Ono equations with arbitrary large initial data in the critical space.} Int. Math. Res. Not. IMRN 2010, no. 2, 297--319.

\bibitem{ZK} V.E. Zakharov and E.A. Kuznetsov,
{\sl On three dimensional solitons.} Sov. Phys. JETP., 39 (1974), 285--286.

\end{thebibliography}
\end{document}